\documentclass[a4paper, 12pt] {article}

\usepackage{fender}  

\begin{document}

    \title{RCF\,3 \\ 
           Map-Code Interpretation via Closure \\
                  $\epsi\,\mathcal{U}\tilde{\eps}$}

\footnotetext{
  this is part 3 
  of a cycle on \emph{Recursive Categorical Foundations}} 

\footnotetext{ 
  Legend of LOGO: \emph{closed} evaluation $\epsi$
  (part of Cartesian Closure) to give, with the help of
  ``stratified'' code \emph{interpretation} into 
  \emph{Universal Chain} $\mcU,$ \emph{code} ``self''-evaluation 
  $\tilde{\eps}.$}

    \author{Michael Pfender\footnote{
              TU Berlin, Mathematik, pfender@math.tu-berlin.de}
           }    

    \date{July 2008\footnote{last revised \today}}                          

\maketitle

\abstract
{For a (minimal) Arithmetical theory with higher Order Objects, 
\ie a (minimal) \emph{Cartesian closed} arithmetical theory -- 
coming as such with the corresponding \emph{closed evaluation} -- 
we \emph{interprete} here \emph{map codes,} out of $\cds{A,B}$ say,
into these maps ``themselves'', coming as elements (``names'') of 
hom-Objects $B^A.$ The interpretation (family) uses a Chain of 
\emph{Universal Objects} $\Un,$ one for each Order stratum with
respect to ``higher'' Order of the Objects. Combined with 
\emph{closed,} axiomatic \emph{evaluation,} these interpretation 
family gives \emph{code-self-evaluation.} Via the usual diagonal 
argument, Antinomie \NAME{Richard} then can be \emph{formalised} 
within our minimal higher Order (Cartesian closed) arithmetical 
theory, and yields this way inconsistency, for all of its extensions, 
in particular of \textbf{set} theories as $\ZF,$ of the Elementary 
Theory of (higher Order) Topoi with Natural Numbers Object as 
considered by \NAME{Freyd} as well as already of the Theory of 
Cartesian Closed Categories with NNO considered by \NAME{Lambek.}
}

\section{Introduction}

Starting point is a discussion of \NAME{Cantor}'s (indirect) 
argument for uncountability of the real numbers 
(in the unit interval), \ie of the \emph{set} $2^\N = \mc{P}\,\N$
of (``actual'' infinit) sequences $a = a(j): \N \to 2.$

This indirect argument assumes all these $a: \N \to 2$ to be
enumerated in form $a_i = a_i(j): \N \to 2,\ i \in \N.$
\NAME{Cantor} then takes as sequence outside this enumeration
of the $a_i$ the sequence 
$\tilde{a} = \tilde{a} (j) \defeq \neg\,a_i(i): \N \to 2.$

But what is this\ $a_i(i) \in 2$\ ? Let us try to apply
\NAME{Cantor}'s argument to any type of \emph{constructive}
real numbers, where in fact there \emph{is} an enumeration,
$a_i$ of all (finite) \emph{texts,} (Computer) \emph{programs,}
standing for -- ``describing'' -- these constructive real numbers,
\eg the primitive recursive power-series descriptions for
$e$ and $\pi.$ But if you want to \emph{change} the diagonal 
\emph{values} in this \NAME{Cantor}'s infinit table $a_i(j)$ 
of the \emph{constructive} reals, you must be able to \emph{evaluate}
the $i$th of these -- say primitive recursive -- programs
at $i \in \N.$ Now \NAME{Ackermann} has shown, that for the case 
of PR function \emph{codes} (``programs'', texts) this 
diagonal evaluation (and then its a posteriori modification) 
cannot be PR any more: The related (equi-complex) 
``Ackermann function'', namely \emph{diagonal evaluation}

\smallskip 
$\eps(f_n,n): \N \xto{\Delta} \N \times \N 
            \xto{\# \times \id_\N} \cds{\N,\N} \times \N \xto{\eps} \N,$
grows faster then any PR function; here $\#(n): \N \xto{\iso} \cds{\N,\N}$
is the PR enumeration of all PR map PR \emph{codes} $f_n$ 
``from'' $\N$ ``to'' $\N.$ The diagonal then says: ``apply''
n-th PR map to -- \emph{evaluate} $n$th PR map code \emph{at} -- 
argument $n.$

$[\,$Presumably this \emph{non-closedness} under code-evaluation 
applies to any constructive class of real numbers and power sets, 
such real numbers obtained e.g by (iterated) ``application'' of 
Intermediate-Value Theorem taken as axiom.$]$

So the possiblity of \emph{closed evaluation,} here of

\smallskip
\quad $\epsi_{\N,2} = \epsi_{\N,2} (\chi,n) = \chi(n) = [\,n \in \chi\,]:
          2^\N \times \N = \mc{P} \times \N \to 2$

\smallskip
is at the basis of classical set theory, with its \emph{closure}
under (iterated) \emph{formation} of power set (and internal
hom sets). This gave rise to investigation of ``all''
the uncountable cardinalities in \text{set} theory, a central 
branch of this theory proper.

The claim of present investigation is that these uncountabilities,
at least a (potentially) infinit ascending chain of uncountabilites, 
leads to a contradiction. The idea is to interpret the \emph{map codes,}
$\code{f} \in \cds{A,B}$ say, of a (minimally presented) theory 
$\PRe$ of PR Arithmetic with (``higher Order'') Cartesian Closure 
added, into these maps ``themselves'', $f \in B^A,$ out
of internal hom Object $B^A,$ in \textbf{set} theory the \emph{map set}

\quad $B^A = \set{f \in \mc{P} (A \times B)\,|
                 \,\forall\,a \in A\ \exists!\,b \in B\ (a,b) \in f}.$

Combined with \emph{closed,} axiomatic evaluation 
$\epsi_{A,B}: B^A \times A \to B,$ $\epsi_{A,B} (f,a) = f(a),$ 
available in \textbf{set} theory and there \emph{needed} for 
(generalisation of) \NAME{Cantor}'s argument above to establish
the strictly ascending hierarchy of \emph{cardinals,} will give a 
\emph{code-self-evaluation,} 
$\tilde{\eps}_{\N,2}: \cds{\N,2} \times \N \to 2,$ and from this
-- because of the \emph{``self''} -- an (anti-)diagonal predicate
$d = d(n): \N \to 2 \xto{\neg} 2,$ whence a \emph{liar} map
$\mathit{liar} = \neg\,\mathit{liar}: \one \to 2$ establishing 
the asserted contradiction for (minimal) Cartesian Closed
PR Theory $\PRe$ and its extensions.     
  
We now outline the sections to come and forshadow at this
occasion some of the notations to be introduced:

  \newpage

\smallskip
\textbf{2\ Theory Closure by Internal hom and Evaluation:}

Here we extend \emph{basic} (categorical) Theory $\PRa = \PR+(\abstr)$ 
of Primitive Recursion with (virtual) extensions $\set{A\,|\,\chi}$
of PR predicates (see part RCF1) by Cartesian Closure, this
in form of adding just new internal hom \emph{Objects,} $B^A,$ 
new \emph{map constants} $\epsi_{A,B}:$ \emph{closed evaluation,} 
and $\axlambda_{A,B}:$ for Cartesian Closure \emph{front adjucntions,} 
as well as suitable \emph{equations} for then already available 
\emph{conjugate} and \emph{coconjugate} maps, but \textbf{no}
new (meta) operations for maps. Resulting Theory is called $\PRe,$
since its decisive ingredient over Theory $\PR$ is closed evaluation
$\epsi_{A,B}: B^A \times A \to B$ with its characteristic equations.

\smallskip
\textbf{3\ Order Stratification for Closed Arithmetic $\PRe$}

In this section we divide higher Order Theory $\PRe$ into
\emph{strata} \\
 $\PRen \bs{\prec} \PRe,$ Cartesian PR theories with 
Order of Objects up to $\ul{n}.$ \textbf{Note:}
$\Ord\,(C^B)^A \defeq \Ord\,C^{B \times A} 
    < \Ord\,C + \Ord\,B + \Ord\,A = \Ord\,C^{B^A} \bydefeq \Ord\,C^{(B^A)},$
``since'' $(C^B)^A \iso C^{B \times A}.$

\smallskip
\textbf{4\ An Ascending, \emph{Universal} Object Chain}

Based on \emph{Universal Object} $\X \subset \N$ for Theory $\PRa,$
$\X$ made out of all (codes of) singletons $\an{n}$ and (possibly nested)
pairs $\an{a;b}$ of natural numbers -- it contains all Objects
$A$ of $\PRa$ coretractively embedded -- we obtain an ascending
\emph{Chain} $\mcU: \mcU_1 = \X \xto{\sqsubset} \mcU_2 = \X^\X \ldots$
of Objects and coretractions, each $\Un$ \emph{universal} for
its stratum $\PRen:$ $A \sqsub \Un$ coretractive for each (pointed)
Object $A$ of Order up to $\ul{n}.$

\smallskip
\textbf{5\ Map-Code Interpretation}

This section develops the central idea of present investigation: An
\emph{interpretation map family}
  $$\intn = [\,\intn_{A,B}: \cds{A,B}_{\PRen} \to B^A\,]_{A,B},
                     \ \ul{n} \in \ul{\N}\ \text{external, ``naive''}$$
is constructed, stratum by stratum, the $\intn$ leading into
Universal Object (at most) $\Unn.$ 

Technically, these Object-pairs indexed families (must and) can 
be ``derived'' from a stratum specific ``global'' Interpretation 
$\Intn = \Intn(u): \Vn \to \Unn,$ $\Vn$ the map code set
of (whole) stratum $\PRen:$ reason for considering 
\emph{Universal Objects,} here: $\Unn.$ 

What we have to do is to \emph{``interprete''} code constants and 
\emph{code operations,} namely (formal) composition, Cartesian 
product and iteration of map codes into the 
\emph{objective correspondants,} \eg -- \emph{plugged} into $\ZF$ --

\smallskip
  $\code{g} \odot \code{f} = \code{g} \code{\circ} \code{f} 
      \bydefeq \code{g \circ f} \overset{\intt} {\mapsto} 
              \intt (g) \circ \intt {f} = g \circ f.$

\smallskip
In our ``formally minimal'' context, this interpretation is based
on the \emph{name} $\name{f}: \one \to B^A$ of a map 
$f: A \to B,$ $\name{f}$ easily defined via \emph{conjugation,}
in \textbf{set} theoretical terms: 
$\name{f} = \set{(\emptyset,f)}: \one \to B^A.$

Interpretation $\intt$ works by the correspondence of 
operations $\odot = \code{\circ},\ \code{\times},$ and
$\code{\S}$ on map \emph{codes} for composition, Cartesian 
product and \emph{iteration} one hand, and associated 
\emph{internal closed} operations, called $\ccode{\circ},$
$\ccode{\times},$ as well as $\ccode{\S}$ on the other. 
These latter are all defined out of set theoretically motivated
``coconjugated'' ones, by \emph{conjugation.} Straightforward
but technically complicated calculations then give the 
central \textbf{Interpretation Theorem,} saying essentially
that (stratum specific) \emph{interpretation}
\begin{align*} 
& \intn_{A,B}: \cds{A,B}_{\PRen} \to B^A \xto{\sqsub} \Unn 
                                 \ \text{is \emph{objective,} \ie:} \\
& \intn_{A,B} (\code{f}) = \name{f}: \one \to B^A,
                             \ \ \text{for}\ f: A \to B \In \PRen.
\end{align*}

\textbf{6\ Self-Evaluation}

With interpretation properties above it is now easy to give a
\emph{sound, objective} code-\emph{self}-evaluation for ``minimal''
Cartesian Closed PR Theory $\PRe,$ namely
\begin{align*}
& \tilde{\eps}_{A,B} (u,a) \defeq \epsi_{A,B} (\intt_{A,B} (u),a): \\
& \cds{A,B}_{\PRe} \times A 
          \xto{\intt \times A} B^A \times A \xto{\in} B,\ \text{with} \\   
& \tilde{\eps}_{A,B} (\code{f},a) = f(a): A \to B. 
                                         & \text{(Objectivity).}
\end{align*}
This then gives immediately formalisation of Antinomie \NAME{Richard}
for $\PRe$ by the usual diagonal argument.

\medskip 
Notions and results for \emph{basic} Theory $\PRa = \PR+(\abstr)$
of Free-Variables (categorical) Theory of Primitive Recursion
with schema of predicate abstraction -- and its Universal Object -- 
are given in \NAME{Pfender}/ \\ 
\NAME{Kr\"oplin}/\NAME{Pape} 1994
and in \NAME{Pfender} 2008 RCF1, RCFX.


\section{Theory Closure by Internal hom and Evaluation} 

We \textbf{extend} here categorical Theory $\PRa = \PR+(\abstr)$ of 
Primitive Recursion -- with \emph{predicate abstraction} 
$\chi \bs{\mapsto} \set{A\,|\,\chi}$ -- into a Theory 
$\PRe \defeq \PRa+(\hom),$ with -- in adddition -- \emph{internal hom} 
$\bfan{A,B} \bs{\mapsto} B^A$ given by \textbf{axiom,} as well as
theory internal -- axiomatic, \emph{closed} -- evaluation    
  $$\epsi = [\,\epsi_{A,B}: 
                    B^A \times A \to B\,]_{A,B\,\in\,\PRe}.$$
This in -- logical -- contrast to \emph{constructive,} 
\NAME{Ackermann} type, formally \emph{partial} -- but still 
\emph{``constructive''} -- evaluation family 
  $$\eps = [\,\eps_{A,B}: \cds{A,B} \times A \parto B\,]_{A,B\,\in\,\PRa}$$
for theories $\pi_O\bfR$ (strengthening $\PRa$) above -- \emph{family} 
obtained out of one single (formally \emph{partial} PR) map
  $$\eps = \eps(u,x): \mrPRa \times \X = \cds{\X,\X}_{\PRa} \times \X 
                                                          \parto \X.$$
\textbf{Comment on Notation:} \emph{Closed evaluation} reads \eg \\  
$\epsi_{A,2} = \epsi_{A,2} (\chi,a) = \chi(a) 
    = [\,a \in \set{A\,|\,\chi}\,]: 
          2^A \times A = \mc{P}A \times A \to 2.$

\smallskip\noindent
This motivates notation for closed evaluation. The ``other'' use
of symbol \text{``\,$\in$''} is -- in Cartesian Theories -- 
\text{``$a \in A\ \free$''}: $a$ is a \emph{(free) variable} on $A,$
categorical meaning: $a$ is (identity of $A$) or a \emph{projection}
onto $A.$ This legitimates free-variables diagram chase below 
categorically. 

Theories $\PRa$ and $\PRe$ fixed, we explain now some 
(known) basic concepts and results, in the language
of Primitive Recursion and Higher Order Arithmetic
sketched above.

\smallskip
Basic for our \emph{Universal Chain} of Objects -- 
\emph{upwards open (!)} -- is the \emph{First Order Universal Object} 
$\X \subset \N$ of all (codes of) \emph{singletons,} $\an{n},$ and 
(possibly nested) \emph{pairs,} $\an{a;b},$ of natural numbers.

Each \emph{fundamental} $\PR$-Object $\one,\N,(\N \times \N)$ \etc 
is \emph{coretractively embedded} into $\X,$ for example \\ 
$(\N \times \N) \owns (m,n) \overset{\sqsubseteq} {\mapsto} 
                    \an{m;n} \in \an{\N \times \N} \subset \X.$  

Extension of Theory $\PRa$ into \emph{Cartesian Closed} Theory 
$\PRe$ \textbf{presented} equationally -- by \NAME{Horn} \ul{inferences} 
-- via additional (formal) \emph{exponential} Objects (Object terms)
of form $(B^A)$ for $A,B$ ``already there'', examples: 
$\N^\N,\ \mcU2 = \X^\X,\ \mcU3 = \X^{(\X^\X)}\ \etc,$ as well
as (additional) families of \emph{map constants} 
\begin{align*}
& \epsi_{A,B}: B^A \times A \to B
    \ (\text{axiomatic, \emph{closed evaluation}}), \\
& \qquad [\,\text{within \textbf{set} theory:}
                \ \epsi(f,a) \bydefeq f(a),\,]\ \text{as well as} \\
& \axlambda_{A,B}: A \to  (A \times B)^B,
    \ \emph{closed front adjunction,}  
      \ \text{``\,$A \owns a \xto{\axlambda} [\,b \mapsto (a,b)\,]$''.}
\end{align*}
These two families are to satisfy the \emph{adjointness} equations for 
(covariant) Functors,
  $A \times B \ \leftadjoint\ B^A: \PRe \bs{\lto} \PRe,
                                       \ (A\ \text{``fixed''}),$
namely defining \emph{conjugation} and \emph{coconjugation} below
as mutually inverse (meta) bijections.

These \NAME{Horn} schemata are merged with those of $\PRa,$ here: 
with forming \emph{Cartesian products} of Objects, with 
\emph{iteration} schema (and Freyd's uniqueness of 
\emph{initialised} iterated), as well as schema $(\abstr)$ of
forming \emph{(virtual) extensions,} cf part RCF\,1. 

\smallskip
Taken together the above internal $\hom$ structure with 
endo map iteration -- and \NAME{Freyd}'s uniqueness of the 
\emph{initialised iterated} -- as well as with (virtual) 
\emph{predicate abstraction} --  
we arrive at Theory $\PRe = \PRa+(\hom) = \PR+(\abstr)+(\hom),$ of 
\emph{Primitive Recursion} with \emph{Object exponentiation} and
\emph{closed evaluation:} Evaluation within the Theory itself. 

\smallskip
$[\,$The latter in contrast to availability of 
``only'' -- \NAME{Ackermann} type, \emph{not} 
PR, (still) \emph{constructive} -- evaluation \textbf{of} 
Theory $\PRa = \PR+(\abstr)$ within ``only'' Theory
$\hatPRa$ --  of formally \emph{partial} PR maps, theory
equivalent to Theory $\muR$ of (partial) $mu$-recursive maps,
see RCF1.$]$


\smallskip
\textbf{Remark:} Theory $\Fin$ of \emph{finite} (number) sets 
\emph{has} internal $\hom$ -- exponentiation -- coming with closed 
evaluation family $\epsi_{A,B}: B^A \times A \to B.$ 
But if you want to \textbf{define} this -- infinitely \ul{indexed} 
family -- made out of (finite) maps, you need Primitive Recursive 
case distinction on $\N \supset B^A,$ and this ``global'',
\emph{mother} evaluation 
  $$\hat{\epsi}: \N \times \N \supset 
            \underset{A,B} \bigoplus (B^A \times A) \to B \subset \N$$
is necessarily \emph{genuine} PR, \emph{not finite.} 

\smallskip
Internal $\hom$ -- and ``closed'' evaluation $\epsi$ -- give, 
within Theory $\PRe$ of \emph{Higher order Arithmetic,} 
\cf \NAME{Eilenberg}\,\&\,\NAME{Kelly} 1966 for internal 
$\hom$ structure, as well as \NAME{Freyd} 1972 and 
\NAME{Lambek}\,\&\,\NAME{Scott} 1986 for the combined structure, 
the following \textbf{defined} map families \emph{conjugation} 
and \emph{coconjugation:} 

\emph{Conjugation} is \textbf{given} by \textbf{schema}  
\inference{ (\emph{conj}) }
{ $f: A \times B \to C$ in $\PRe$ }
{ $\ol{f} = \mr{conj} [f] \defeq f^B \circ \axlambda_{A,B}: 
                               A \to (A \times B)^B \to C^B,$ \\
& in \textbf{set} theory conjugate $\ol{f}$ realised  as \\ 
& $a \overset{\ol{f}} \mapsto [b \mapsto (a,b) \mapsto f(a,b) \in C\,],$
}
and \emph{coconjugation} is introduced by \textbf{schema}
\inference{ (\emph{coconj})}
{ $g: A \to C^B$ in $\PRe$}
{ $\ol{g} = \mr{coconj} [g] \defeq \epsi_{B,C} \circ (g \times \id_B):
                                  A \times B \to C^B \times B \to C,$ \\
& in \textbf{set} theory coconjugate $\ol{g}$ realised as \\
& $[\,(a,b) \mapsto (g(a),b) \mapsto \ol{g}\,(a,b) 
                                 \bydefeq (g(a)) (b) \in C\,].$
}

These two families are to satisfy -- by \textbf{axiom,} and do so
(already within \emph{finite}) \textbf{set} theory and the 
\emph{Elementary Theory of Topoi} \textbf{ETT}
in place of Theory $\PRe$ around to be ``constructed'' -- the following 
\emph{higher order} meta-\ul{bi}j\ul{ection} \ul{e}q\ul{uations}: 
\inference{ (\emph{co/conj})}
{ $f: A \times B \to C$ in $\PRe$}
{ $\PRe \derives\ \ol{\ol{f}} = \mr{coconj} [\mr{conj} [f]] = f: 
                                                A \times B \to C$}
as well as
\inference{ (\emph{conj/co})}
{ $g: A \to C^B$ in $\PRe$}
{ $\PRe \derives\ \ol{\ol{g}} = \mr{conj} [\mr{coconj} [g]] = g: 
                                                         A \to C^B.$}


\smallskip
The above data, in particular (axiomatically given) \emph{families} 
$\axlambda$ and $\epsi,$ \textbf{define} the following 
meta-\ul{ma}p, and make it into a \emph{covariant functor} $\hom$
-- the \emph{covariant internal hom functor} -- via the following
schema:

\smallskip
\cinference{ (\hom\text{-co}) }
{ $A,\ g: B \to C$ in $\PRe$ }
{ \cinference{} 
  { $g \circ \epsi: B^A \times A \xto{\epsi} B \xto{g} C$ }
  { $g^A \defeq \ol{g \circ \epsi_{A,B}}: B^A \to C^A$ }
}

\medskip
Analogeous schema \textbf{defining} the \emph{contravariant} (closed) 
internal hom \emph{functor:}

\medskip
\cinference{ (\hom\text{-contra}) }
{ $A,\ g: B \to C$ in $\PRe$ }
{ \cinference{} 
  { $\epsi\,\circ\,(A^C \times g): 
      A^C \times B \xto{A^C \times g} A^C \times C \xto{\epsi} A$ }
  { $A^g \defeq \ol{\epsi\,\circ\,(A^C \times g)}: A^C \to A^B.$ }
}

\medskip
All four: \emph{Universal} property, the two Functor properties, 
and \emph{right adjointness,} of covariant closed internal hom
 $g \bs{\mapsto} g^A$ -- namely right adjointness to 
\emph{Cylindrification} 
  $$\bfan{g: B \to C} \bs{\mapsto} 
      \bfan{A \times g: A \times B \to A \times C},
                        \quad\text{Object $A$ fixed},$$
are consequences of the pair $\mr{conj}/\mr{coconj}$ above 
to be a pair of meta-\ul{bi}j\ul{ections}, \ul{inverse} to 
each other.

\smallskip
\textbf{Remark;} $\axlambda_{A,B}$ and $\epsi_{A,B}$ are natural 
transformations, but we will not rely on these properties here.


\section{Order Stratification for Closed Arithmetic $\PRe$} 

\textbf{Definition:} The -- formal -- Order $\ul{\mr{Ord}}\,A$ of a 
\emph{higher order Object} -- of Theory $\PRe$ -- is \textbf{defined} 
externally PR as follows:
\begin{align*}
& \Ord\,\one,\ \Ord\,\N \defeq 1, \\
& \Ord\,(A \times B) \defeq \max\set{\Ord\,A,\Ord\,B}, \\
& \Ord\,\set{A\,|\,\chi: A \to 2} \defeq \Ord\,A, \\
& \qquad 
    \text{in particular}\ \Ord\,2 = \Ord\,\set{n \in \N\,|\,n < 2} = 1, \\
& \qquad
    \Ord\,\X = \Ord\,\set{\N\,|\,\X: \N \to 2} = 1 
                      \ (\X\ \text{is a predicative subset of}\ \N.) \\
& \text{For}\ B \In \PRa\ \text{and}\ A \In \PRe\ 
                                (\Ord\,A\ \text{``already known''}): \\
& \Ord\,B^A = 1+\Ord\,A; \\
& \quad \text{finally: for}\ C \in \PRa,\ B,C \In \PRe: \\
& \Ord\,(C^B)^A \defeq \Ord\,C^{B \times A} 
                                 \bydefeq 1+\Ord\,(B \times A) \\
& \bydefeq 1+\max(\Ord\,B,\Ord\,A).
\end{align*}
The latter clause takes in account the (canonical) $\PRe$ 
\emph{reduction} isomorphism $(C^B)^A \iso C^{B \times A}.$ 

With this \textbf{definition,} we have in particular
$\Ord\,B^A \leq \Ord\,B+\Ord\,A$ for all $\PRe$ Objects $A,B,$
as well as $\Ord\,\mcU\ul{n} = \ul{n},$ \eg
$\Ord\,\mcU3 \bydefeq \Ord\,\X^{\X^\X} = \Ord\,\X^{(\X^\X)} = 3.$

So subSystem $\PRa$ of Theory $\PRe$ has all its (presenting) Objects 
of Order 1, it is our basic, ``1st'' Order, subSystem of Theory 
$\PRe$ -- not a priori an (``embedded'') sub\emph{Category,} since 
the higher-order axioms of $\PRe$ may entail -- within $\PRe$ -- 
\emph{new} equations between map terms of $\PRa$ viewed as map terms 
of $\PRe,$ in \emph{logical} terms: The \emph{Extension} $\PRe$ of 
$\PRa$ may be not \emph{conservative.}

\smallskip
\textbf{Broadening} to \textbf{Theories Extension Chain:} 
We \textbf{define} an exhaustive Chain of \ul{subS}y\ul{stems}
  $\PRen\,\bs{\preceq}\,\PRe,\ \ul{n} \in \ul{\N},$ 
\ PR as follows:

\smallskip\noindent
-- $\PRe\,1 \defeq \PRa;$ \\
-- Assume $\PRen\,\bs{\preceq}\,\PRe$ 
to be known via its (canonical) \emph{presentation:} 

\smallskip
\emph{Object} terms, \emph{map} terms, schemata 
for \emph{map} (term) \emph{equations.}

\smallskip
Then subSystem $\PRen+ = \PRe[\ul{n}+1]$ is 
\textbf{defined} to be the \emph{Cartesian-PR-Category Closure} 
of subSystem $\PRen$ \ul{mer}g\ul{ed} with
Closure under \emph{formal adjunction} of 
\begin{itemize}
\item
all Objects of Order $\ul{n}+1$ 

\item
the canonical isomorphisms $(C^B)^A \xto{\iso} C^{B \times A}$ 
given in $\PRe$ for $C \In \PRa,$ $A,B \In \PRen,$ and their 
inverses $C^{B \times A} \xto{\iso} (C^B)^A$
 
\item
$\PRe$ \ul{families}
  $\axlambda_{A,B}: A \to (A \times B)^B$ as well as 
     $\epsi_{A,B}: B^A \times A \to B,$ \\
this for $\Ord\,A + \Ord\,B,\ 2\,\Ord\,B \leq \ul{n}+1,$ and
$\Ord\,A + \Ord\,B \leq \ul{n}+1$ respectively.


\end{itemize}

Additional (merged) \textbf{equations} come in, for the maps of 
$\PRen+,$ via schemata $(\mr{co/conj})$ as well as $(\mr{conj/co})$ 
of $\PRe$ (above), which are to establish the 
\emph{conjugation/coconjugation} \ul{bi}j\ul{ection} for all those 
of their instances, for which all formal ingredients -- Object terms 
and map terms -- are \ul{enumerated} so far within $\PResn.$ 



\medskip
\textbf{Corollary} to this \textbf{Definition:} 
\begin{enumerate} [(i)]

\item \emph{Conjugation upgrade:}
\inference{ (\mr{upgrade}) }
{ $f: A \times B \to C$ in $\PRen\,\bs{\prec}\,\PRe,$ }
{ $\ol{f} = \mr{conj} [f] = f^B \circ \axlambda_{A,B}:  
                                    A \to (A \times B)^B \to C^B$ \\ 
& lives in $\PRenn \bs{\prec} \PRe.$ 
}

\item \emph{Coconjugation upgrade:} 
\inference{ (\text{co-upgrade}) }
{ $g: A \to C^B$ in $\PRen\,\bs{\prec}\,\PRe$ }
{ $\ol{g} = \mr{coconj}\,(g) = \epsi_{B,C} \circ (g \times \id_B):$ \\
& $A \times B \to C^B \times B \to C$\ \  lives already in $\PRen:$
}

Critical exponential Object $C^B$ is presupposed to 
belong already to Theory $\PRen.$ 

\item
Theory $\PRen$ contains Objects up to Order $\ul{n},$
and in fact \ul{some} of its Objects \emph{have} this Order.

\item 
External \ul{ascendin}g\ul{ ``union''} of \ul{all} subSystems 
$\PRen,\ \ul{n} \in \ul{\N},$ \emph{exhausts} Theory
$\PRe,$ \ie gives a -- \emph{stratified} -- \textbf{presentation} 
of Theory $\PRe:$ Objects, maps, and equations.


\end{enumerate}


\section{An Ascending, \emph{Universal} Object Chain} 

Basic -- 1st Order -- Arithmetical Theory $\PRe\,1 = \PRa$ 
\emph{has} a \emph{Universal Object} in itself, a 
\emph{first-Order Universal Object,} namely 
the Object $\X \subset \N$ -- of (codes of) all singleton 
(lists) and of pairs, possibly nested: binary bracketed NNO 
tuples.

$\X$ is a Universal Object -- of Theory $\PRa$ and therefore also of its 
\emph{stengthenings,} as for example for the full first order
\emph{sub}category $\barPRe\,1$ of $\PRe.$ Object $\X$ is 
\emph{universal} in the following sense:

$\X$ admits -- for each $\PRa$-Object $A,$ an \emph{embedding} 
(here an \emph{injective map}), even a \emph{coretractive} map
(see below), $\sqsub_A\,: A \xto{\sqsub} \X,$ \textbf{defined} 
externally PR in the obvious way.

All these embeddings 
 $\sqsub_A\,: A \ovs{\sqsub} \X$ 
-- \emph{disjoint} as far as \emph{fundamental} Objects $A$ 
are concerned, namely binary bracketed powers of $\N,$ no
genuine \emph{abstracted} sets -- come with canonical 
\emph{retractions}
  $\sqsup_A: \X \ovs{\sqsup} A,$
the latter equally for abstracted Objects $\set{A\,|\,\chi}$
having a \emph{point,} $a_0: \one \to \set{A\,|\,\chi},$
as in particular $\X \subset \N,$ coming with ``its'' \emph{zero}
$\an{0}: \one \to \X.$ 

\medskip
\textbf{Graded-Universal-Object Chain:} Each of our Theories
$\PRen$ in the hierarchy -- except (!) ``roof'' Theory 
$\PRe$ itself -- comes with a \emph{canonical} 
Universal Object, $\Un = \upexp{\X} {\ul{n}},$ 
externally PR \textbf{defined} as follows, as an internal
version of a \NAME{Grothendieck}-Universe (?):
\begin{align*}
& \U_1 = \upexp{\X} {1} \defeq \X^1 = \X, \\
& \Usn \defeq \X^{\Un} = \X^{\upexp{\X} {\ul{n}}} 
                               \bydefeq \upexp{\X} {\ul{n}+1}  
\end{align*}
For opening the possibility that a \emph{higher, later} 
Universal Object in the chain is good also as 
Universal Object for a \emph{lower, earlier} Theory in the 
hierarchy, we establish first the \emph{Universal Chain} 
$\mc{U}$ as a chain of \emph{embeddings} 
$\sqsub\ =\ \sqsub_{\ul{n}}: \Un \to \Usn$
coming each with a retraction \\ 
$\sqsup\ =\ \sqsup_{\ul{n}}\,: \Usn \to \Un,$ as follows:

\smallskip
Universal Chain $\mc{U}$ begins with (commutative) \textsc{diagram}
$$
\xymatrix{
& & \X^{\one}
    \ar[dll]_{\epsi_{\X}}^{\iso}
\\
\U_1 = \X
\ar @<-0.5ex> [rrrr]_{\sqsub}
\ar[drr]_{\bar\ell_{\X,\one}}^{\iso}
& & & & \X^\X = \U_2
        \ar @<-0.5ex> [llll]_{\sqsup}
        \ar [ull]_{\X^{\one \ovs{\an{0}} \X}} 
\\
& & \X^{\one}
    \ar [urr]_{\X^{\X \ovs{!} \one}}      
}
$$
Diagram chase in case of set theory:
$$
\xymatrix{
& & {[\,0 \mapsto \an{0} \mapsto x\,]}
    \ar @<-0.7ex> @{|->} [dll]_{\epsi_{\X}}^{\iso}
\\
x
\ar @<-0.5ex> @{|->} [rrrr]_{\sqsub}
\ar @<-1.0ex> @{|->} [drr]_{\bar\ell_{\X,\one}}^{\iso}
& & & & {[\,y \mapsto 0 \mapsto x\,]}
        \ar @<-0.5ex> @{|->} [llll]_{\sqsup}
        \ar @<-0.7ex> @{|->} [ull]_{\X^{\one \ovs{\an{0}} \X}} 
\\
& & {[\,0 \mapsto x\,]}
    \ar @<-1.0ex> @{|->} [urr]_{\X^{\X \ovs{!} \one}}      
}
$$

\smallskip
The general Universal Chain member then is recursively 
\textbf{defined} by commutativity of \textsc{diagram}

$$
\xymatrix{
\Un = \X^{\U_{\ul{n}-1}}
\ar @<-0.5ex> [rrrr]_{\sqsub}
\ar @/_2pc/ [rrrr]_{\X^{\Un\,\sqsup\,\U_{\ul{n}-1}}}  
& & & & \X^{\Un} = \Usn
        \ar @<-0.5ex> [llll]_{\sqsup}
        \ar @/_2pc/ [llll]_{\X^{\U_{\ul{n}-1}\,\sqsub\,\Un}} 
}
$$
Easy Diagram chase for verifying section/retraction property \eg in
\textbf{set} theory.


\smallskip
Generalising the above to the case of $B^A$ instead of $\X^{\U_{\ul{n}-1}}$
we now \textbf{define} recursively the (coreteractive) embeddings 
  $$\sqsub\ =\ \sqsub_{B^A}\,: B^A \to \Usn,\ B^A
                               \In \PResn = \PRe\,[\ul{n}+1],$$
based on the (coretractive) embeddings $\sqsub_B\,: B \into \X = \U_1$
above, as follows, ``but'' first only for Object $B$ in $\PRe1 = \PRa:$

\smallskip
-- Anchor: for $A$ in $\PRa,$ (natural)
embedding $\sqsub_A\,: A \into \mcU_1 \bydefeq \X^1 = \X$
has been \textbf{defined} above by converting natural
numbers $n$ in singleton codes $\an{n},$ and -- recursively --
pairs in code pairs, out of $\PRa$ \emph{Universal Object}
$\X \subset \N.$ Furthermore, a canonical \emph{retraction}
$\sqsup_A\,: \X \to A$ for the embedding has  been mentioned above, 
for Object $A$ coming with a \emph{point,} $a_0: \one \to A$ say.

\smallskip
-- Step: Assume embedding $\sqsub_A\,: A \into \Un$ to be given, 
together with retraction $\sqsup_A\,: \Un \onto A,$ for ``each'' 
Object $A$ of Order $\ul{n}$ -- in $\PRen.$

Consider then a (genuine) Object in $\PResn,$ of form $B^A,$ 
$A$ in $\PRen,$ $B$ in $\PRa$ (!). Then the  \textsc{diagram} 
below -- simplified one of the former one above -- \textbf{defines} 
``universal'' embedding and retraction for Object $B^A$ 
into/from $\Usn \bydefeq \X^{\Un}:$
$$
\xymatrix{
& & B^{\Un}
    \ar[dll]_{B^{\sqsup_A}} 
\\ 
B^A
\ar @<-0.5ex> [rrrr]_{\sqsub_{B^A}}
\ar[drr]_{{\sqsub_B}^A}
& & & & \X^{\Un} = \Usn
        \ar @<-0.5ex> [llll]_{\sqsup_{B^A}}
        \ar [ull]_{{\sqsup_B}^{\Un}}
\\
& & \X^A
    \ar[urr]_{\X^{\Un\,\sqsup\,A}}
}
$$
Again easy Diagram chase for verifying section/retraction property
in case of \textbf{set} theory.

The general, not \emph{normal form} case, of a $\PResn$ Object of 
form $B^A,$ $B = D^C$ not basic, not in $\PRa,$ is reduced to the 
above one via (natural) isomorphism $(D^C)^A \iso D^{(C \times A)}$ 
-- such isomorphism possibly applied several times --, to a 
\emph{normal form} case Object to be embedded, by a map \emph{within}
$\PResn$ (or lower) \emph{into} $\Usn$ or lower, by the method above 
for the case of Object $B$ in $\PRa.$ Embedding into $\Usn$ in the 
latter case then is by composition with embedding 
$\mcU_{\ul{m}} \xto{\sqsub} \Usn,$ $\ul{m} < \ul{n}+1.$

Taken together the above -- including the modification for the
non-normal-form case -- we have $\PResn$ embedded all Objects 
of $\PResn$ into $\Usn,$ namely all $\PRe$ Objects of -- up to -- 
Order $\ul{n}+1.$ This \textbf{proves}

\bigskip
\textbf{Embedding Theorem} for Chain $\mc{U}:$ 
\begin{enumerate} [(i)]

\item Each single of our Theories $\PRen$ 
admits \emph{coretractive} embeddings 
    $\sqsub_A: A \ovs{\sqsub} \Un$
for each of its (pointed) Objects $A,$ into ``its'' 
\emph{Universal Object} within the \emph{section/retraction} \emph{Chain}
$$
\xymatrix{
\mc{U}: \quad 
\mcU_1 = \X\phantom{i}
\ar @<-0.5ex> @{>->} [r]
& \quad
  \ar @<-0.5ex> @{->>} [l]
  \ar @<-0.5ex> @{>->} [r]
  & \Un\phantom{i}
    \ar @<-0.5ex> @{->>} [l]
    \ar @<-0.5ex> @{>->} [r]
    & \Usn\phantom{i}
      \ar @<-0.5ex> @{->>} [l]
      \ar @<-0.5ex> @{>->} [r]
      & {}
        \ar @<-0.5ex> @{->>} [l]
}     
$$      
of these ``Universal'' Objects, the Chain $\mc{U}$ hosted 
as an ascending chain in \emph{global,} \ul{hi}g\ul{her Order} 
Theory $\PRe.$





\item By the above discussion of -- canonical -- natural 
retractions $\sqsup_{\ul{n}}\,: \Usn \to \Un,$ retractions 
to embeddings $\sqsub_{\ul{n}}\,: \Un \to \Usn,$ the above 
coretractive embedding for all Objects of $\PRen,$ 
into $\Un,$ gives also (canonical) embeddings into 
\emph{later} Objects of chain $\mc{U},$ \ie if $\Un$ is 
replaced by $\mcU_{\ul{m}},\ \ul{m} > \ul{n},$ and (coretractive) 
\emph{embedding} $A \to \mcU_{\ul{m}}$ is taken as
\ $\sqsub_A: A \ovs{\sqsub} \Un \ovs{\sqsub} \ldots 
                                 \ovs{\sqsub} \mcU_{\ul{m}}.$
\end{enumerate}



\section{Map-Code Interpretation} 
 
Using Order Stratification above -- of 
\emph{higher order Cartesian Closed Theory} $\PRe = \PRa+(\hom)$ --
we now \textbf{define} -- via PR -- a Theory-internal
\emph{interpretation} map family

\smallskip
  $\intn = [\,\intn_{A,B}: \mapn{A}{B} 
              \defeq \cds{A,B}_{\PRen} \to B^A\,]_{A,B},
                                         \ \ul{n} \in \ul{\N},$
 
\smallskip
$A,B$ Objects of \emph{stratum} $\PRen;$ interpretation 
$\intn_{A,B}$ will be \textbf{defined} inside stratum $\PRenn.$ 

\smallskip
\textbf{Example:} 
  $\intt^1_{\N,2}: \cds{\N,2}_{\PRa} = \cds{\N,2}_{\PRe\,1} \to 2^\N$
will live inside stratum $\PRe2$ -- and higher --, see discussion
in foregoing section.

\smallskip
$[\,$Such a \emph{stratum} is a PR Cartesian Theory, but it is 
\emph{truncated} what concerns (exponential) Order 
of Objects and (axiomatic) evaluation. We will see below -- 
in particular for our \emph{interpretation} of \emph{constructive,} 
PR defined ``internal'' hom sets $\cds{A,B}$ into \emph{closed}
ones $B^A,$ that it is sufficient to climb up to stratum $2\,\ul{n}$
for interpretation of stratum $\ul{n}.]$

\smallskip
In our present -- categorical -- context, $\ulfamily$ 
$\intt_{A,B} = \intn_{A,B},$ $\ul{n} \in \ul{\N}$ 
\ul{fixed}, \emph{can} and \emph{must} (?) be 
\textbf{defined} formally as (a family) derived from \emph{one} 
single $\PRenn$ map. 
So, as one Interpretation for all -- on stratum
$\PRen$ \ul{fixed} -- we are lead to \textbf{define} 
-- PR over $\PRenn$ -- this \emph{global} Interpretation 
as a $\PRenn$ map, with suitable, \emph{universal,} 
\emph{Domain} and \emph{CoDomain.} 
 
We start by \emph{type-description} of this 
$\ulfamily$ -- to be defined, later, as a $\ulfamily$ of 
\emph{Domain/Codomain} restrictions of the 
\emph{one single} map $\Int^{\ul{n}}$ of Theory 
$\PRen$ to be (objectively) PR \textbf{defined} -- of 
following type:
 
  $\Int^{\ul{n}} = \Int^{\ul{n}} (u):
     \Vn \defeq \underset{A,B} {\overset{\quad} {\bigoplus}}\,\mapn{A}{B}
                  \xymatrix{{} \ar @{-->} [r] & \Unn,}$ where
$\mapn{A}{B}$ is an abbreviation for internal, 
\emph{syntactical} $\PRen$-\emph{map code (!) set} 
  $\cds{A,B}_{\PRen} \subset \mr{V} \subset \N,$ \emph{from} $A$ to $B,$
$A,B$ both Objects of $\PRen.$


\smallskip
We turn now our ``typifying'' proposal (!) above, into 
a \textsc{diagram} which displays a special -- central -- countable 
sum (``disjoint union''), and its (litteral) component-\emph{inclusions.} 
This ``special'' \emph{sum}-\textsc{diagram} is available within 
$\PRen$ -- as \emph{litteral,} disjoint union of 
\emph{predicates,} disjoint by definition. 

\smallskip
\emph{Global} $\PRen$ \emph{Interpretation}
$\Intn,\ (\ul{n}\ \text{fixed}),$ to be \textbf{defined} following 
actual type-discussion, then will be \textbf{characterised} a 
posteriori (!) as $\PRenn$ map, \emph{induced} map out of the 
(countable) sum, induced by its \emph{components} 
$\intn_{A,B}: \mapn{A}{B} \to B^A \sqsubset \Unn,$
$A,B$ in stratum $\PRen.$ 

In other words: $\Intn$ will be \emph{PR ``constructed''} -- ``over'' 
$\PRen,$ ``but only'' within $\PRenn$ -- in such a way 
that it becomes the (unique) $\PRenn$ map out of \textbf{sum} 
$\Vn \subset \N,$ which makes \emph{commute} the following 
(externally) \ul{countable} \textsc{diagram,} this diagram 
\emph{available} within $\PRenn:$ 

\bigskip
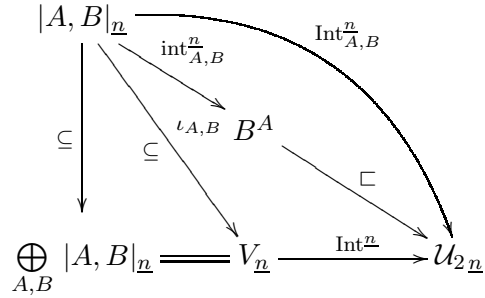
\begin{minipage} {\textwidth}
\xymatrix{
\mapn{A}{B}
\ar @/^3pc/ [rrrdd]^{\Intn_{A,B}}
\ar [dd]_{\subseteq}
\ar [rdd]^{\iota_{A,B}}_{\subseteq}
\ar [rd]^{\intn_{A,B}}
& &
\\
& B^A
  \ar [rrd]^{\sqsubset}
  &
\\
\overset{\phantom{M}} {\underset{A,B} {\bigoplus}}\ \mapn{A}{B}
\ar @{=} [r]
& \Vn
  \ar [rr]^{\Intn} 
  & & \Unn & 
}

\bigskip           
Interpretation map \textsc{diagram} $(A,B \In \PRen)$
\end{minipage}

\bigskip


PR Construction of $\PRen$ map 
$\Intn: \N \supset \Vn = \underset{A,B} {\bigoplus}\,\cds{A,B} \to \Unn$ 
is recursively \emph{merged} with that of maps 
  $\intn_{A,B}: \mapn{A}{B} \to B^A,$ 
the latter being (recursively) \textbf{defined} as 
Domain/Codomain restrictions of \emph{universal} PR 
defined Interpretation map $\Intn$ within $\PRenn,$ 
in fact by the followig \textbf{defining} commutative \textsc{diagram}
($B$ pointed):
$$
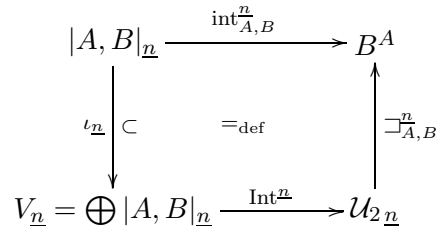
\xymatrix @+2em{
\mapn{A}{B}
\ar @{} [rd] |{\defeq}
\ar[r]^{\intn_{A,B}}
\ar[d]_{\iota_{\ul{n}}}^{\subset}
& B^A
\\
\Vn = \bigoplus \mapn{A}{B}
\ar[r]^(0.6){\Intn}
& \Unn
  \ar[u]_{\sqsup^{\ul{n}}_{A,B}} 
}
$$
This type of restriction becomes \emph{possible} -- at least 
easier -- by the fact that ``all'' maps considered come as
\emph{section/retraction} pairs. This is in particular the case
for all \emph{injections-into-sums} embeddings here to be treated.



\smallskip
\textbf{Constructive Internalisation} of meta operations for
our Theories $\PRe$ and subSystems $\PRen:$
 
\emph{Composition} $\circ: \T \bs{\times} \T \bs{\lto} \T(A,C)$ 
of $\T$ -- Theory $\T$ any (categorical) theory -- 
\emph{constructively} internalises to 

  $\odot = \code{\circ}: \cds{B,C}_\T \times \cds{A,B}_\T \to \cds{A,C}_\T,
               \ (v,u) \mapsto\ \an{v \odot u} \in \cds{A,C}_\T.$
As Objects $A,B$ here all Objects of $\T$ are allowed, for 
$\T :\,= \PRen$ in particular Object $\Un$ and its (embedded) subobjects.
 
Analogeously Cartesian product ``$\times$'' has as \emph{coded version} 
family 
  $\cds{A,B}_\T \times \cds{C,D}_\T 
      \owns (u,v) \overset{\code{\times}} {\mapsto} 
        \an{u \code{\times} v} 
          \in \cds{A \times C,B \times D}_\T,$
for (arbitrary) $\T$-Objects $A,B,C,D,$ including in particular
Objects $\Un$ in case of theory $\T :\,= \PRen.$



Analogeously for \emph{iteration} ``$\S$'' within (Cartesian) PR theories 
in particular ``again'' for extension $\PRen$ of PR Theory 
$\PRa = \PR+(\abstr):$ 
  $\cds{A,A} \owns v \overset{\code{\S}} {\mapsto} v^{\code{\S}} 
                                        \in \cds{A \times \N,A},$
here \eg for iteration of $\PRen$ endo maps with Domain 
$A :\,= \Un$ and their internalisations. 

\smallskip
\textbf{Definition:}
The \emph{constructive $\PRen$-codes} in $\Vn$ 
are -- \textbf{first} -- the \emph{constructive} internal 
\emph{map-constants} 

\smallskip
  $\code{\axlambda_{A,B}}: 
                \one \to \mapn{A}{(A \times B)^B}, 
                                          \ \text{and}\ 
       \code{\epsi_{A,B}}:  
                  \one \to \mapn{B^A \times A}{B}$

\smallskip
for $\axlambda_{A,B},\ \epsi_{A,B}$ in $\PRen.$

\smallskip
\textbf{Second:} the ``derived'' \emph{Cartesian map constants} for 
the new Objects and their Cartesian products -- with the ``old'' ones 
and with the new ones --: identities, terminal maps, (left and right) 
projections, and 

\textbf{Third:} ``Closure'' under \emph{composition} and 
\emph{cylindrification} (Cartesian product with an identity)  
as well as under \emph{iteration} of endo maps.

 
\smallskip   
Next we \textbf{define,} for $f: A \to B$ in $\PRen$ 
-- and hence in particular Objects $A,B$ in $\PRen,$ 
the notion \emph{name of} $f: A \to B,$ symbolised as 
  $\name{f} = \name{f: A \to B}: \one \to B^A,$
available in stratum $2\,\ul{n}.$ 

This up-to-2\,\ul{n}th Order construct 
$\name{f}$ is defined simply by \emph{conjugation,} as
  $\name{f} = \ol{f \circ r_{\one,A}} 
      = \mr{conj} 
          [\,f \circ r: \one \times A 
             \xto{\iso} A \to B\,]: \one \to B^A.$
 
\emph{Name} $\name{f}$ of $f$ represents, meta-bijectively, 
map $f: A \to B$ \emph{within} -- as \emph{defined element} 
of -- \emph{closed} internal $\hom$ set $B^A.$
In \textbf{set} theory: 
  $\name{f} \defeq \set{(\emptyset,f)}: 
          \one \to B^A \subset \mc{P} (A \times B).$ 

\smallskip
By its \textbf{definition} via \emph{conjugation,} $\name{f}$
has characteristic property 
$\epsi_{A,B} (\name{f},a) = f(a) = f: A \to B,$ $a \in A\ \free.$
Verification of this \emph{closed Objectivity} from definition is 
trivial for set theoretic environment, and straight forward for 
the general higher Order case. 

\smallskip
\textbf{Definition} of \emph{global} Interpretation
  $\Intn = \Intn(u): \N \supset \Vn \to \Unn,$ \\
of (\ul{n}-truncated), internal \emph{map-code-set} 
$\Vn$ of Theory $\PRen$ into $\PRenn$'s  
\emph{Universal Object} $\Unn$ -- within (the language of) Theory 
$\PRenn,$ is by \emph{recursive case distinction} on the 
\emph{structure} of the \emph{map code} $u \in V_{\ul{n}}$ 
to be \emph{interpreted.} (At beginning we do not typify into 
types $A,B$ for $\mapn{A}{B} \subset \Vn.$)
  %
This PR \textbf{case distinction} for Definition
of \emph{Interpretation} 
  $\Intn(u): \Vn \to \Unn$ runs as follows:

\smallskip
-- \textbf{Case} of $\PRen$ \emph{map constants} ``$\bas$'', namely
$0: \one \to \N$ and $s: \N \to \N$ as well as \emph{all} Cartesian map 
constants of $\PRen:$ identities, terminal maps, diagonals, (binary) 
projections, as well as case of the additional -- \emph{closed} -- 
map constants $\axlambda_{A,B},\ \epsi_{A,B}$ of $\PRen:$

\smallskip
For all of these \emph{anchor cases,} we \textbf{define} 
\emph{Interpretation} $\Intn = \Intn(u): \Vn \to \Unn$ 
in the below -- PR -- by ``codes to names:'' 
\begin{align*}
& \Int(\code{\bas}) \defeq \name{\bas}: \one \to \Unn,\ \eg \\
& \Int(\code{\ell: A \times B \to A}) 
           \defeq \name{\ell: A \times B \to A}: 
               \one \to A^{A \times B} \ovs{\sqsub} \Unn, \\
& \text{Objects}\ A,B\ \text{in stratum}\ \PRen.
\end{align*}
This gives in particular for the ``extra'' $\PRen$ basic codes, 
with appropriate $\PRen$ Objects as \emph{types:} 
\begin{align*}
& \Int(\code{\axlambda_{A,B}: A \to (A \times B)^B})
                              \defeq \name{\axlambda_{A,B}}: \\ 
& \qquad
    \one \to ((A \times B)^B)^A \ovs{\sqsub} \Unn,\ \text{as well as} \\
& \Int(\code{\epsi_{A,B}: B^A \times A \to B})
    \defeq \name{\epsi_{A,B}}: 
              \one \to B^{B^A \times B} \ovs{\sqsub} \Unn,
\end{align*} 
The latter two ``inclusions'' $\ \sqsub \Unn$ are available
by the fact that $\axlambda_{A,B}$ and $\epsi_{A,B}$ were supposed 
to live ``already'' within $\PRen,$ and that conjugation -- at
the base of \emph{name} $\name{f}$ -- at most doubles Order of 
(minimal) ``receiving'' stratum, here Order $\ul{n}.$ 



\smallskip
What we still have to worry about is \emph{self-referential (!)} 
Interpretation of family members $\intn_{A,B}: \cds{A,B} \to B^A,$ 
obtained from \\
$\Intn: \Vn = \mr{PR}\epsi\,\ul{n} \to \Unn$ 
by \emph{Domain/CoDomain restriction.}

For these \emph{injections} into \emph{sum} $\Vn,$ we \emph{will} 
obtain (!), out of our PR \textbf{case-}definition of \emph{global} 
Interpretation $\Intn: \Vn \to \Unn,$ by 
\textbf{definition} -- below -- of \ul{families} 
\begin{align*}
& \intn_{A,B}: \mapn{A}{B} \to B^A
         \ \text{as Domain/CoDomain \emph{restrictions} of}\ \Intn: \\
& \Intn({\code{\intn_{A,B}}}) = \name{\intn_{A,B}}: 
        \one \to (B^A)^{\cds{A,B}} \iso B^{A \times \cds{A,B}} 
                                           \ \text{within}\ \PRenn.
\end{align*}
The latter map will lead in fact -- Order verification -- into 
$\PRenn$ by our \textbf{definition} of 
\begin{align*}
& \Ord\,(B^A)^{\cds{A,B}} = \Ord\,B^{A \times \cds{A,B}} \\
& \leq \Ord\,B+\max(\Ord\,A,\,\Ord\,\cds{A,B}) \\ 
& \qquad = \Ord\,B+\max(\Ord\,A,\,\Ord\,\N) \\
& \leq \Ord\,B+\Ord\,A \leq 2\,\ul{n},
\end{align*}
and since the isomorphism pair 
$(B^A)^{\cds{A,B}} \iso B^{A \times \cds{A,B}}$
is included in $\PRenn$ by \textbf{definition} of stratum 
$\PRen+ = \PRe[\ul{n}+1].$ 





 
\smallskip      
Based on the anchor cases above, we define by 
\emph{genuine primitive recursion} stratum Interpretation
$\Intn$ of (constructively) \emph{composed} codes, 
\emph{Cartesian ``parallelised''} as well as of \emph{iterated} 
ones, as follows by PR \textbf{case} distinction on 
\emph{Iteration Domain} for PR \textbf{definition} of 
$\Intn: \Vn \to \Unn,$ PR case distinction on 
the \emph{disjoint components} $\mapn{A}{B}$ of 
``syntactic (code) universe'' $\Vn \subset \N,$ 
which in turn is a PR defined predicative subObject of 
$\N$ within Theory $\PRa$ -- in the r\^ole of (internal)
\emph{\ul{Metamathematics}} -- $\PRa$ \ul{subS}y\ul{stem} 
of $\PRen\,\bs{\prec}\,\PRe\ [\ = \text{``\,$\PRe\ul{\infty}$\,''}\ ].$


\smallskip
With -- always below -- \textbf{abbreviation}

\smallskip
  $\mapn{A}{B} \bydefeq \cds{A,B}_{\PRen} \subset \Vn = \mrPRen \subset \N,$
we introduce 

\smallskip
$\PRen$ map (map-$\ulfamily$, indexed on $\ul{n} \in \ul{\N}$) 
  $$\Intn = \Intn(u): 
    \Vn = \underset{A,B} {\overset{\quad} {\bigoplus}}\,\mapn{A}{B} 
                                                           \to \Unn,$$
\emph{merged} with its \emph{Domain/Codomain} restrictions, recursively 
as follows:

\smallskip
Interpretation of \emph{constructive} internal \emph{composition:}
For $A,B,C$ in stratum $\PRen:$
\begin{align*}
& \text{for}\ u \in \mapn{A}{B} \subset \Vn,
               \ v \in \mapn{B}{C} \subset \Vn \\ 
& \qquad\ [ \implies \an{v \odot u} \in \mapn{A}{C} \subset \Vn\ ]: \\ 
& \Intn \an{v \odot u} 
    \quad [\ = \intn_{A,C} (\an{v \odot u}) \subset \Vn\ ] \\ 
& \defeq \Intn(v) \ccode{\circ} \Intn(u) \\   
& \bydefeq \ccode{\circ}\,(\Intn(v),\Intn(u)): \\ 
& \Vn \times \Vn \xto{\supset} \mapn{B}{C} \times \mapn{A}{B} 
                \xto{\intn_{B,C} \times \intn{A,B}} C^B \times B^A \\
& \xto{\ccode{\circ}} C^A \xto{\sqsub} \Unn.
\end{align*}
This is a formally defined $\PRenn$ map, in particular since
$\Vn \times \Vn \xto{\supset} \mapn{B}{C} \times \mapn{A}{B}$
is -- obviously -- a \emph{retraction.} 
We recall further that ``embedding'' $C^A \xto{\sqsub_{C^A}} \Unn$
also comes with a retraction, $\Unn \xto{\sqsup_{C^A}} C^A.$

\smallskip
\emph{Axiomatic} internal composition -- \emph{competing} 
with \emph{constructive} internal composition 
$\odot = \code{\circ},$ gets a similar symbol, $\ccode{\circ},$ 
which may be read \emph{Closed internal composition,} similarly: 
$\ccode{\times}:$ \emph{Closed internal Cartesian product,} 
as well as $\ccode{\S}$ for \emph{Closed internal iteration:}
For the general background on \emph{Closed Categories} see 
\NAME{Eilenberg}\,\&\,\NAME{Kelly} 1966.

\smallskip
\emph{Closed internal composition} 
  $\ccode{\circ}: \Unn \xto{\sqsup} C^B \times B^A 
                            \xto{\ccode{\circ}_{A,B,C}} C^A,$ \\
(\emph{retraction} 
     $\sqsup\ =\ \sqsup_{C^B \times B^A}$ just cares on -- feasable --
case distinction) 
 
\smallskip
is \textbf{defined} via \emph{conjugate} 
  $\ccode{\circ} = \ol{\ol{\ccode{\circ}}}$ of
\begin{align*}
& \ol{\ccode{\circ}}_{A,B,C}: (C^B \times B^A) \times A \to C \\
& \defeq [\,\epsi\ \circ\,(C^B \times \epsi):
              C^B \times B^A \times A 
                \xto{C^B \times \epsi} C^B \times B 
                                                \xto{\epsi} C\,], 
\end{align*}
with Cartesian \emph{associativity} (natural) isomorphisms of form 
  $$A \times B \times C \defeq (A \times B) \times C 
          \underset{\iso} {\xto{\mathrm{ass}}} 
                                      A \times (B \times C)
                                           \ \text{omitted.}$$
This \textbf{case} of $\PRe$-map 
$\Intn: \Vn \xto{\supset} \mapn{A}{C} \to C^A \xto{\sqsub} \Unn$
describes in fact a $\PRenn$ map:
 
In its chain of Objects -- and in its Order \ul{minimal}
presentation of maps -- it is at most of Order 2\,\ul{n} -- for 
Objects $A,B,C$ all of Order at most $\ul{n}.$

\smallskip
-- \emph{Interpretation} of constructive internal 
\emph{product of maps:}
This is analogeous to the above, even easier, since 
the two components of a Cartesian product are completely independent
of each other, ``exercise''.

\bigskip
 
-- \textbf{Case} of an internally \emph{iterated} 
$v^{\code{\S}} \in \mapn{A \times \N} {\N},$ $v \in \mapn{A}{A}$ free, 
Object $A$ in $\PRen.$
\textbf{Define} in this case 
\begin{align*}
& \Vn \supset \mapn{A}{A} \owns v \overset{\Intn} {\mapsto} 
          \Intn(v^{\code{\S}}) \in A^{A \times \N} \xto{\sqsub} \Unn 
                                                         \ \text{by} \\
& \Intn(v^\S) \defeq \ccode{\S} (\Intn(v)): 
     \Vn \xto{\supset} \mapn{A}{A} 
           \xto{\ccode{\S}} A^{A \times \N} \xto{\sqsub} \Unn.
\end{align*}  
Here $\PRen{2}$ map $\ccode{\S}: A^A \to A^{A \times \N}$ is defined
as \emph{conjugate} to 
\begin{align*}
& \ol{\ccode{\S}} = \ol{\ccode{\S}}_A: A^A \times (A \times \N) \to A, \\ 
& \qquad \text{this in turn \textbf{defined} -- PR -- by} \\
& \ol{\ccode{\S}}_A (v,(a,0)) \defeq a: 
                         A^A \times (A \times \one) \to A, \\
& \ol{\ccode{\S}}_A (v,(a,s\,n)) 
    \defeq \epsi_{A,A}\,(v,\ol{\ccode{\S}} (v,(a,n))): 
                              A^A \times (A \times \N) \to A.
\end{align*}  


With the above, in particular with definition of \emph{Interpretation}
map $\Intn$ on \emph{map constants} -- among them (the codes of) 
$\axlambda$ and $\epsi,$ $\Intn$ is (PR) defined on all of its 
arguments, in particular on \emph{conjugated} and $\hom$-functor
values, since these are definable in terms of Composition,
Cartesian Product and Iteration out of the 
\emph{basics.} Furthermore, the above \emph{type insertions} show 
that $\PRen$ map 
  $$\Intn = \Intn(u): \N \supset \Vn 
                = \bigoplus\,\mapn{A}{B} \to \Unn,
                                 \ \Ord{A} \leq \ul{n},$$ 
is -- as expected -- induced by \emph{Object-pair typified family} 
  $$\intn_{A,B}: \mapn{A}{B} \to B^A \xto{\sqsub^{\ul{n}}_A} \Unn,$$
($\ul{n}$ still \emph{fixed}), more precisely: it is the 
\emph{induced} out of \ul{countable} sum:
\begin{align*}
& \Intn = (\,\sqsubset_{B^A} \circ\, \intn_{A,B}: 
                  \mapn{A}{B} \to B^A \sqsub \Unn\,)_{A,B}: \\  
& \qquad \xymatrix{\bigoplus\,\mapn{A}{B} \ar @{-->} [r] & \Unn.}
\end{align*} 
By \textbf{Definition} of \emph{constructive coding} -- namely by 
definition of code \emph{composition} 
$v \odot u = v \code{\circ} u,$ 
of code \emph{product} $u \code{\times} v,$ 
and of code \emph{iteration} 
$u^{\code{\S}},$ all simply given by concatenation of 
\ul{ASCII} strings -- we have the following

\medskip
\textbf{Structure Preservation} by \emph{Constructive Coding:}
\begin{align*}
& \qquad \emph{Composition:}\ \text{for}\ f: A \to B
                 \ \text{and}\ g: B \to C \In \PRen: \\ 
& \code{g\,\circ\,f} = \code{g} \odot \code{f} 
               \bydefeq \code{g} \code{\circ} \code{f}: \\ 
& \one \to \mapn{A}{C} \bydefeq \cds{A,C}_{\PRen}; \\
& \qquad \emph{Cartesian product:}\ \text{for} \ f: A \to C
                 \ \text{and}\  g: B \to D \In \PRen: \\
& \code{(f \times g)} = \an{\code{f} \code{\times} \code{g}}:
                \one \to \mapn{A \times B} {C \times D}, \\
& \qquad \text{as well as}\ \emph{Iteration:}\ \text{for}
            \ f: A \to A \In \PRen: \\
& \code{f^\S} = \code{f}^{\code{\S}}: \one \to \mapn{A \times \N} {A}.
\end{align*}
For \emph{closed} internalisation we have an analogeous result, namely

\smallskip
\textbf{Structure Preservation} by \emph{Closed Internalisation:}
\emph{Naming} 
  $$\bfan{f: A \to B} \bs{\mapsto} \bfan{\name{f}: \one \to B^A}$$
\emph{preserves} Composition, map-Product and iteration
\emph{into} the corresponding closed families 
  $\ccode{\circ}_{A,B,C},\ \code{\times}_{A,B,C,D},
                       \ \text{as well as}\ \code{\S}_A,$
in detail: 

\smallskip
-- \emph{Composition:} For $A \xto{f} B \xto{g} C$ in $\PRen$ we have:
\begin{align*}
& \name{g \circ f} = \ccode{g} \ccode{\circ} \ccode{f} 
      \bydefeq \ccode{\circ}\,(\name{g},\name{f}): \\ 
& \quad \one \xto{(\name{g},\name{f})} C^B \times B^A 
                                    \xto{\ccode{\circ}} C^A,
\end{align*}
it lives within stratum $\PRenn.$

\smallskip
-- \emph{Cartesian product:} 
For $f: A \to C,$ and $g: B \to D$ in $\PRen:$
\begin{align*}
& \name{(f \times g)} = \an{\name{f} \ccode{\times} \name{g}}
                           \defeq \ccode{\times} (\name{f},\name{g}): \\ 
& \quad \one \to C^A \times D^B 
                   \xto{\ccode{\times}} (C \times D)^{A \times B},
\end{align*}
this again lives in stratum $\PRenn.$

\smallskip
-- \emph{Iteration:} For $f: A \to A$ in $\PRen,$
  $$\name{f^\S} = \ccode{\S} (\name{f}) = \ccode{\S} \circ \name{f}: 
                                    \one \to A^A \to A^{A \times \N},$$
it is likewise a  $\PRenn$ map.

\smallskip
\textbf{Proof:} 

-- (Central), \emph{Composition case:} 
We consider first \emph{coconjugated} composition, namely
\begin{align*}
& \ol{\ccode{\circ}} \circ\,((\name{g},\name{f}) \times A): 
        \one \times A \to C^B \times B^A \times A 
                                   \xto{\ol{\ccode{\circ}}} C \\
& \bydefeq \epsi \circ\,(C^B \times \epsi) 
                       \circ\,((\name{g},\name{f}) \times A): \\
& \qquad \one \times A \to C^B \times B^A \times A 
               \xto{C^B \times \epsi} C^B \times B \xto{\epsi}  C \\
& = g \circ f \circ r: \one \times A \to B \to C. & (*) 
\end{align*}
The latter equation follows from the evaluation properties of 
\emph{closed evaluation} instances $\epsi: B^A \times A \to B,$ 
and $\epsi: C^B \times B \to C,$ 
by \emph{Free Variable} chasing -- namely 
free variable $a :\,= r_{\one,A}: \one \times A \onto A.$ 

\smallskip
By \emph{conjugation} of (both sides of) the above equation we get the 
assertion in the present composition case:
\inference{ (\emph{conj})}
{ $\ol{\ccode{\circ}}\,\circ\,((\name{g},\name{f}) \times A):
                     \one \times A \to C^B \times B^A \times A \to C$ \\
& $= (g \circ f)\,\circ r: \one \times A \to A \to C$ 
}
{ $\ccode{\circ}\,\circ\,(\name{g},\name{f}): 
                             \one \to C^B \times B^A \to C^A$ \\
& $= \name{g\,\circ\,f}: \one \to C^A.$
}

\smallskip
-- Case of \emph{Cartesian product:} analogeous, ``exercise''. 

\smallskip
-- \emph{Iteration case:} We start again with the 
\emph{conjugate} side: For a $\PRen$ endo $f: A \to A,$ 
we want to show
\begin{align*}
& \ol{\name{f^\S}}: \one \times (A \times \N) \to A \\
& \bydefeq f^\S \,\circ\,\iso\,: 
       \one \times (A \times \N) \xto{\iso} A \times \N \xto{f^\S} A \\
& = \ol{\ccode{\S}}\,\circ\,(\name{f} \times (A \times \N)):   & (***) \\
& \one \times (A \times \N) \to A^A \times (A \times \N) 
                                        \xto{\ol{\ccode{\S}}} A.
\end{align*} 
For \textbf{Proof} of $(***)$ we use the \textbf{definition} above, of
$\ol{\ccode{\S}},$ \emph{coconjugate} of 
$\ccode{\S}: A^A \to  A^{A \times \N},$
making commute the lower two rectangles of the following diagram:
$$
\xymatrix@1{
& & A \times \N
    \ar @/^1pc/ [rrd]^{f^\S} 
    \ar @{} [d] |{(***)} \\                               
\one \times (A \times \N) 
  \ar @/^1pc/ [urr]^{\iso}
  \ar[rr]^-{\name{f} \times (A \times \N)}
  \ar[d]^{\iso}
  \ar @{} [rrd] |{=}
  & & A^A \times (A \times \N)
      \ar[rr]^-{\ol{\ccode{\S}}}
      \ar[d]^{\iso}
      \ar @{} [rrd] |{\bydefeq}
      & & A  \\
(\one \times A) \times \N
  \ar[rr]^-{(\name{f} \times A) \times \N}
  & & (A^A \times A) \times \N
      \ar[rr]^-{(\ell,\epsi)^\S}
      & & A^A \times A
          \ar @{>>} [u]_{r}
}
$$
For showing $(***),$ we show commutativity of the frame diagram, by free
variables diagram chasing, with free variables $a :\,= \ell_{A,\N}$ 
$n :\,= r_{A,\N}:$  

$$
\xymatrix@1{
& & (a,n)
    \ar @/^1pc/ @{|->} [rrd]^{f^\S} 
                              \\                               
(0,(a,n)) 
  \ar @/^1pc/ @{|->} [urr]^{\iso}
  \ar @{|->} [d]^{\iso}
  & & 
      & & f^\S(a,n)  \\
((0,a),n)
  \ar @{|->} [rr]^-{(\name{f} \times A) \times \N}
  & & ((\name{f},a),n)
      \ar @{|->} [rr]^-{(\ell,\epsi)^\S}_{(\mapsto)}
      & & (\name{f},f^\S(a,n))
          \ar @{|->} [u]_{r}
}
$$
Remains to show $(\mapsto),$ \ie to show:
\begin{align*}
& (\ell,\epsi)^\S ((\name{f},a),n) = (\name{f},f^\S(a,n)): & (\bullet) \\ 
& (\one \times A) \times \N \to A^A \times A. 
\end{align*}
We show this by \ul{external} Peano Induction, \ie by 
\emph{uniqueness of the iterated,} as follows:
\begin{align*}
& (\ell,\epsi)^\S ((\name{f},a),0) 
     = (\name{f},a) = (\name{f},f^\S(a,0))   & (\anchor) \\
& \text{as well as} \\
& (\ell,\epsi)^\S ((\name{f},a),n+1) 
     = (\ell,\epsi)^\S ((\name{f},\epsi_{A,A}\,(\name{f},a)),n) \\
& = (\ell,\epsi)^\S ((\name{f},f(a)),n) \\
& \qquad \text{by \emph{ evaluation property} of}\ 
                               \epsi_{A,A}: A^A \times A \to A \\
& = (\name{f},f^\S(f(a),n))\ \text{by induction hypothesis on}\ n \\
& = (\name{f},f^\S(a,n+1)): A \times \N \to A.  & (\step)
\end{align*}
This shows $(\bullet),$ \ie $(\mapsto)$ in the \textsc{diagram:} 
Map $(\ell,e)^\S$ -- diagram -- throws in fact $((\name{f},a),n)$ 
into $(\name{f},f^\S(a,n)).$ So assertion $(***)$ above has been 
\textbf{shown}. Whence, by \emph{conjugation:}
\begin{align*}
& \name{f^\S} = \mr{coconj} [\,f^\S \,\circ\,\iso\,: 
      \one \times (A \times \N) \to A \times \N \xto{f^\S} A\,] \\
& = \ccode{\S}\,\circ\,\name{f}:
      \one \xto{\name{f}} A^A \xto{\ccode{\S}} A^{A \times \N}, 
\end{align*}
and that \textbf{proves} the remaining case of 
\emph{Structure Preservation} via Closed Internalisation\ \,\textbf{\qed}  


\smallskip
We now come to our central result, the

\bigskip
\textbf{Interpretation Theorem:}
\begin{enumerate} [(i)]

\item
\emph{CoDomain Suitability} of interpretation $\ulfamily$:
PR defined $\PRenn$ \emph{interpretation family} 
$\intn_{A,B}: \mapn{A}{B} \to \Unn$ -- indexed by Object-pairs, 
\emph{stratum} (strata) $\PRen$ (and $\PRenn$) -- restricts 
in its (single) \emph{CoDomains} to 
\begin{align*}
& \intn_{A,B}: \mapn{A}{B} \to B^A \ [ \ \xto{\sqsub} \Unn\ ] & (*)
\end{align*}
within $\PRenn,$ in form of a commuting \textsc{diagram,} for
$B$ having a point:

\bigskip
\begin{minipage} {\textwidth} 
\xymatrix @+2em{ 
& \mapn{A}{B}
  \ar[ld]_{\subset} 
  \ar[r]^{\intn_{A,B}}
  \ar @/^1pc/ [d]^{\subset}
  & B^A
    \ar @/^1pc/ [d]^{\sqsub_{B^A}}
\\
\underset{A,B} {\overset{\quad} {\bigoplus}} \mapn{A}{B}
\ar @{=} [r]
& \Vn
  \ar @/^1pc/[u]^{\supset}
  \ar[r]_{\Intn}
  & \Unn
    \ar @/^1pc/ [u]^{\sqsup_{B^A}}
}

\bigskip
Interpretation \textsc{diagram:} stratum by stratum, \\ 
\quad 
  global/individual with respect to map-code sets 
\end{minipage}

\bigskip

\item
\emph{Objectivity} within one \emph{stratum:} For $f: A \to B$ in 
$\PRen \bs{\subseteq} \PRe,$ we have
\begin{align*}
  \PRen \derives\ 
  & \intn_{A,B}\,(\code{f}) \bydefeq \intn_{A,B}\,\circ\,\code{f} 
                                      = \name{f}: \one \to B^A & (**)
\end{align*}
\emph{Codes ``originating from'' Objective level are interpreted into names.}

\item
\emph{Stratum-Globalisation} of Interpretation: 
Stratum-indexed $\ulfamily$ 

  $[\,\intn_{A,B}: \mapn{A}{B} \to B^A\,]_{\ul{n} \in \N}$
admits, within Theory 
 
  $\PRe =\ \,\underset{\ul{n}} {\overset{\quad} {\upunion}} \PRen,$ 
Object $\cds{A,B}$ of $\PRe$ again as an \emph{ascending Union,} written
  $\cds{A,B} \bydefeq \cds{A,B}_{\PRe} 
      =\ \, \underset{\ul{n}} {\overset{\quad} {\upunion}}\mapn{A}{B},$
predicatively, and has the universal property of an 
\emph{inductive limit} by PR ``construction''. 

In particular, $\ulfamily$ $\intn_{A,B}: \mapn{A}{B} \to B^A$
above induces a -- \emph{unique} -- \emph{strata-global} map
  $\intt_{A,B}: \cds{A,B} \to B^A\ (***)$ 
\ making commute the following \textsc{diagram:}

\bigskip
\begin{minipage} {\textwidth}
$
\xymatrix{
{\ldots}
\ar[r]^{\subset}
& \mapn{A}{B}
  \ar[rr]^{\subset}
  \ar[rd]^{\subset}_{\iota_{\ul{n}}}
  \ar @/_1pc/ [rddd]_{\intn_{A,B}}
  & & \cds{A,B}_{\ul{n}+}
      \ar[ld]_{\iota_{\ul{n}+}}^{\subset}
      \ar @/^1pc/ [lddd]^{\intt^{\ul{n}+}_{A,B}}
      \ar[r]^{\subset}
      & {\ldots}
\\
& & \qquad \cds{A,B}\,=\ \underset{\ul{n} \in \ul{\N}} {\upunion} \mapn{A}{B} 
    \ar @{-->} [dd]^{\intt_{A,B}}
\\
\\
& & B^A
}
$
\begin{center} Strata-global interpretation \textsc{diagram} \end{center}
\end{minipage}

\bigskip

\item
Strata-global \emph{Objectivity} of Interpretation, ``Codes to names'':

For an arbitrary $\PRen$ map $f: A \to B$ we have:
\begin{align*}
\PRe \derives\ 
& \intt_{A,B}\,(\code{f}) \bydefeq \intt_{A,B}\,\circ\,\code{f} \\ 
& = \name{f}: \one \to B^A \xto{\sqsub} \Unn.  & (\bullet)
\end{align*}

\end{enumerate}

\textbf{Proof:}
\begin{enumerate} [(i)]
\item Type control 
  $\Vn \supset \mapn{A}{B} \owns u 
            \mapsto \intn_{A,B} (u) \in B^A \sqsub \Unn:$

\smallskip
This is \textbf{proved} by structural induction on $u,$ \ie
on $\depth(u): \cds{A,B} \supset \mapn{A}{B},$ $\ul{n}$ ``suitable''
such that all the finitely many \emph{building blocks} $v,w,\ldots$
are in finitely many components of \emph{sum} \\
  $\Vn = \underset{A,B} {\overset{\quad} {\bigoplus}} \mapn{A}{B}.$
This type assertion has been (pre-) discussed already above.

\item
\textbf{Proof} of second assertion $(**)$ 
on \emph{Objectivity} of each member of the $\ul{n} \in \ul{\N}$
and Object-pair $A,B$ indexed family is now as expected, namely
by external \ul{structural induction} on (external) 
$\uldepth[f]$ of map $f: A \to B$ in $\PRe$ in question,
with $\code{f} \in \mapn{A}{B},$ suitable $\ul{n} \in \ul{\N}:$
Each such $f$ comes with such a ``suitable'' $\ul{n},$ since
obviously the $\PRen,\ \ul{n} \in \ul{\N},$ exhaust all
of Theory $\PRe$ here considered.

Now here is the \textbf{Proof} of \emph{Interpretation-Objectivity,}
by structural induction on $\uldepth[f: A \to B]$ ``to be
interpreted'':

\smallskip
For $f$ one of the \emph{map constants} of $\PRe = \PRa+(\hom)$
with $\uldepth[f] = 1$ say, in particular for the members of 
adjunction map families $\epsi_{A,B}$ and $\axlambda_{A,B},$ 
the assertion is trivial, by definition of interpretation 
$\Int,$ and corresponding $\intn_{A,B}$ in these cases.
We now consider $\PRe$ maps with greater $\uldepth:$


\smallskip
For\ $f: A \to B,\ g: B \to C$ in $\PRe$ 
\begin{align*}
\PRen \derives\ 
& \intn_{A,C}\,(\code{g\,\circ\,f}) 
      = \intn_{A,C}\,\an{\code{g} \odot \code{f}} \\
& \qquad [\,\odot\, = \code{\circ}\ 
                 \text{\emph{constructively} internalises}\ \circ\,\,] \\
& \bydefeq \intn_{B,C}\,(\code{g}) \ccode{\circ} \intn_{A,B}\,(\code{f}) \\
& = \name{g} \ccode{\circ} \name{g}\ 
                         \text{by \ul{h}yp\ul{othesis} on $f$ and $g$} \\
& = \name{g\,\circ\,f}: \one \to C^A,
\end{align*} 
the latter by the \emph{composition case} of 
\emph{Structure preservation by Axiomatic Internalisation} 
above.

\smallskip
Similar (\ul{recursive}) \textbf{Proof} for the assertion in 
case of the other binary meta-\ul{o}p\ul{eration}, the 
Cartesian product of maps.

\smallskip
-- Remains the case of an \emph{iterated} 
$f^{\S}: A \times \N \to A,$ given by the unary 
meta-\ul{o}p\ul{eration} $\S$: 
In this case we have     
\begin{align*}
\PRen \derives\ 
& \intn_{A \times \N,A}\,(\code{f^{\S}})  
       = \intn_{A \times \N,A}\,(\code{f}^{\code{\S}}) \\
& \qquad \text{by definition of constructive code of an iterated} \\
& \bydefeq \ccode{\S} (\intn_{A,A}\,(\code{f})) \\  
& \qquad \text{(``homomorphic'' PR definition 
                         of interpretation $\Int$)} \\
& = \ccode{\S} \,\circ\,\name{f}\ \text{by hypothesis on}\ \uldepth [f] \\
& = \name{f^{\S}}: \one \to A^{A \times \N},
\end{align*}
the latter, eventually, by the \emph{iteration case} of 
\emph{Structure Preservation of Closed Internalisation.}  
  
\end{enumerate}
The last two \textbf{assertions} of the \textbf{Theorem}
 -- $(***)$ and $(\bullet)$ -- 
follow straightforward from the former two, by the 
inductive-limit property of our Universal Chain $\mc{U}.$ 

\smallskip
\textbf{Comment:} The members of family
$\Intn: \Vn \to \Unn$ are special maps -- Objective 
map terms -- of Theory $\PRenn \bs{\prec} \PRe,$ and are therefore
covered ``themselves'' by the -- in this regard 
\emph{self-related} \textbf{Interpretation Theorem} above.
This is the reason why I have choosen
as a \emph{Universal Class} \emph{not} a single Object 
or ``super-Object'' for Theory $\PRe,$ but an 
\emph{ascending chain} of ``Universal Objects'' $\Un,$ 
such that Object $\Unn,$ of $\PRenn,$ hosts in particular
\emph{interpretation} of all map codes of stratum $\PRen:$
Chain $\mcU$ is ``upwards open'', think at \NAME{Hilbert}'s hotel. 


\section{Self-Evaluation}     

Here is the key \textbf{Consequence} of the two last assertions
$(***)$ and $(\bullet)$ of the \textbf{Interpretation Theorem,} namely 
possibility for a \emph{constructive self-evaluation} of Theory $\PRe:$ 

\smallskip
\textbf{Define} code-\emph{self}-evaluation family 
\emph{for} Theory $\PRe,$ called  
  $\tilde{\eps} = \tilde{\eps}_{A,B}: \cds{A,B}_{\PRe} \times A \to B,$ 
\emph{within} Theory $\PRe$ itself as  
\begin{align*}
& \tilde{\eps}_{A,B} = \tilde{\eps}_{A,B}\,(u,a) 
                 \defeq \epsi_{A,B}\,(\intt_{A,B}\,(u),a): \\  
& V \times A \supset \cds{A,B}_{\PRe} \times A 
               \xto{\intt_{A,B} \times A} B^A \times A \xto{\epsi} B.
\end{align*}
\textbf{Comment:} Here we used assertion $(***)$ for availability 
of suitable Order-global interpretation family 
  $$\intt_{A,B}: \cds{A,B}_{\PRe} =
     \ \,\underset{\ul{n}} {\overset{\quad} {\upunion}} \mapn{A}{B} 
                                                         \to B^A.$$
We get further, by last assertion -- $(\bullet)$ -- of the 
Theorem, \emph{objectivity} of self-evaluation $\tilde{\eps},$ namely:
for (any) $f: A \to B$ in $\PRe$
\begin{align*}
\PRe \derives\ 
& \tilde{\eps}_{A,B}\,(\code{f},a) = \epsi_{A,B}\,(\intt_{A,B}\,(\code{f}),a) \\
& = \epsi_{A,B}\,(\name{f},a) = f(a): A \to B. & (\bs{*})
\end{align*}
For this latter equation see introduction -- and discussion -- of 
\emph{name} of $f$ above,
  $\name{f} = \name{f: A \to B}: \one \to B^A$ -- in \textbf{set} theory:
$\name{f} = \set{(\emptyset,f)}: \one \to B^A \subset \mc{P} (A \times B).$


\smallskip
Based on this \emph{self-evaluation} family 
of Theory $\PRe,$ we now find within $\PRe$ the following 
(anti) diagonal $d = d(n): \N \to 2:$ 
$\PRe$-map $d = d(n): \N \to 2$ is \textbf{defined} as

  \quad$d \defeq \neg\ \circ\ \tilde{\eps}_{\N,2}\,\,\circ\,(\#,\id_\N): 
      \N \lto \cds{\N,2}_{\PRe} \times \N 
                \xto{\tilde{\eps}} 2 \xto{\neg} 2,$

with $\# = \#(n): \N \xto{\iso} \cds{\N,2}_{\PRe}$ the 
-- isomorphic -- PR \emph{count} of all (internal) predicate codes
(``Klassenzeichen'' in G\"odel's sense), of Theory $\PRe.$ 
As expected in such \emph{diagonal argument} -- Antinomie Richard quoted
by \NAME{G\"odel} -- we \emph{substitute,} within Theory $\PRe,$ the 
counting index 
  $q \defeq \#^{-1}\,(\code{d}) = \#^{-1} \circ \code{d}: 
                \one \to \cds{\N,2}_{\PRe} \xto{\iso} \N,$
of $d$'s \emph{code} into $\PRe$-map $d: \N \to 2$ itself, 
and get a ``liar'' map $\liar: \one \to 2,$ called \emph{liar} 
because it turns out that this map is its own negation, as follows:
\begin{align*}
\PRe\ \derives\ 
& \liar \defeq d\,\circ\,q: \one \to \N \to 2 \\
& \bydefeq d\,\circ\,\#^{-1}\,\circ \code{d} \\
& \bydefeq \neg\,\circ\,\tilde{\eps}_{\N,2}\,\circ\,(\#\,,\,\id_\N) 
               \,\circ\,\#^{-1}\,\circ \code{d} \\
& = \neg\,\circ\,\tilde{\eps}_{\N,2}
               \,\circ\,(\code{d}\,,\,\#^{-1}\,\circ \code{d}) \\
& \bydefeq \neg\,\circ\,\tilde{\eps}_{\N,2}\,(\code{d}\,,\,q) \\
& = \neg\,\circ\,d(q) = \neg\,\circ\,d\,\circ\,q & (\bs{**}) \\
& \bydefeq \neg\ \liar: \one \to 2 \to 2,
\end{align*}
a contradiction: The argument is equation marked $(\bs{**}),$
which is a special instance of \emph{objectivity} equation $(\bs{*})$ 
above, objectivity of \emph{self-evaluation} $\tilde{\eps},$ 
which has been \textbf{defined} within theory $\PRe$ out of
\emph{closed evaluation} $\epsi$ composed with interpretation
family $\intt,$ of map \emph{codes} into \emph{names.}

\medskip
\textbf{Conclusion:}
The argument shows incompatibility of (even just \emph{potential}) 
\emph{infinity} with (formally, axiomatically given) 
\emph{Cartesian Closed} ``Higher Order'' structure of Theory $\PRe.$ 


We obtain this way inconsistency of all extensions of 
Theory $\PRe,$ in particular of -- higher order -- \textbf{set} 
theories, and also of any type of higher Order Arithmetic,
even when given in a categorical setting, as in particular
in \NAME{Lawvere} 1963, and then in \NAME{Freyd}'s 1972 setting 
of (higher Order) Topos Theory with NNO, and in that of 
\NAME{Lambek}\,\&\,\NAME{Scott} 1986. 

The present argument does not depend on \emph{quantification} nor
on availability of a subobject classifier: the (equality) 
\emph{predicates} we rely on here are given by the Cartesian PR 
Arithmetic of theories considered.

\smallskip
\textbf{Disclaimer:}
The argument does \emph{not} apply to \emph{Closed Categories} 
in the sense of \NAME{Eilenberg}\,\&\,\NAME{Kelly,} since there 
is no NNO required for the theory. In the applications, \eg 
Categories of Modules, there is an NNO only \emph{downstairs,} 
in a suitably conceived category of sets.
But that NNO does not bear (naturally) the structure of an
abelian group. 

Even if you consider the category of abelian \emph{semi}-groups
which includes semigroup $\N = \bfan{\N,0,+}:$ an \emph{iterated}
$f^\S: A \times \N \to A$ will not become linear, even not 
bilinear, and hence even not linear when converted into 
a map $f^\S: A \tensor \N \to A$ from the tensor product into $A.$
So this category cannot have $\N$ as an NNO in any suitable way.

Analogeously, the original Elementary Theory $\ETT$ of Topoi
seems me to be not concerned, $\ETT$ in the sense explained by 
\NAME{Wraith} 1973 on the base of mainly (?) \NAME{Lawvere} 1970, 1972, 
and \NAME{Tierney} 1971, as well as more recently explained in
\NAME{Lawvere} \& \NAME{Shanuel} 1991: 

The data and axioms for this \emph{genuine} Theory of Topoi do not 
include an NNO. The motivating examples for Topoi are Categories 
of \emph{sheaves} over a topological space. 
Question: Do these -- Cartesian Closed -- Categories 
come with an NNO on sheaf level? By the above, they cannot come so, 
except they are based on an -- inconsistent -- Cartesian Closed 
set Theory with NNO.



\medskip
\textbf{Problem:}
Diagonal map above is a map within subSystem $\PRen,$
subSystem of Theory $\PRe,$ for $\ul{n}$ from \ul{some} 
$\ul{n}_0$ upwards. Presumably an upper bound for such 
contradictory Order $\ul{n}_0$ can be \ul{calculated}. 
It would be certainly interesting to know a lower bound
$\ul{n}_0$ making $\mc{P}^{\ul{n}}$ contradictory, incompatible
with (potential) \emph{infinity,} in the sense of availability of 
a Natural Numbers Object $\N.$

\section* {References}



  
  
  



 

  \NAME{S.\ Eilenberg, C.\ C.\ Elgot} 1970: \emph{Recursiveness.}
  Academic Press.
  
  \NAME{S.\ Eilenberg, G.\ M.\ Kelly} 1966: Closed Categories. 
  \emph{Proc.\ Conf.\ on Categorical Algebra}, La Jolla 1965, pp. 421-562. 
  Springer. 
  
  \NAME{S.\ Eilenberg, S.\ Mac Lane} 1945: General Theory of Natural 
  Equivalences. \emph{Trans.\  AMS} 58, 231-294.
  
  
  \NAME{P.\ J.\ Freyd} 1972: Aspects of Topoi. 
  \emph{Bull.\ Australian Math.\ Soc.} \textbf{7,} 1-76.
  
  \NAME{K.\ G\"odel} 1931: \"Uber formal unentscheidbare S\"atze der
  Principia Mathematica und verwandter Systeme I. 
  \emph{Monatsh.\ der Mathematik und Physik} 38, 173-198.
  

  \NAME{R.\ L.\ Goodstein} 1971: \emph{Development of Mathematical
  Logic,} ch. 7: Free-Variable Arithmetics.  Logos Press.
  



  \NAME{D.\ Hilbert}: Mathematische Probleme. Vortrag Paris 1900. 
  \emph{Gesammelte Abhandlungen.} 
  Springer 
  1970. 
  
  
  \NAME{P.\ T.\ Johnstone} 1977: \emph{Topos Theory.} Academic Press

  \NAME{A.\ Joyal} 1973: Arithmetical Universes. Talk at Oberwolfach.
  
  
  \NAME{J.\ Lambek, P.\ J.\ Scott} 1986: \emph{Introduction to higher order 
  categorical logic.} Cambridge University Press.
  
  
  \NAME{F.\ W.\ Lawvere} 1964: An Elementary Theory of the Category of
  Sets. \emph{Proc.\ Nat.\ Acad.\ Sc.\ USA} \textbf{51,} 1506-1510.

  \NAME{F.\ W.\ Lawvere} 1970: Quantifiers and Sheaves. 
  \emph{Actes du Congr\`es International des Math\'ematiciens.}
  Nice, pp.\ 329-334.

  \NAME{F.\ W.\ Lawvere, S.\ H.\ Shanuel} 1997 (1991): 
  \emph{Conceptual Mathematics, A first introduction to categories.}
  Cambridge University Press.

  \NAME{S.\ Mac Lane} 1972: \emph{Categories for the working mathematician}. 
  Springer.
  
  
  \NAME{B.\ Pareigis} 1969: \emph{Kategorien und Funktoren}. Teubner.
  
  \NAME{R.\ P\'eter} 1967: \emph{Recursive Functions}. Academic Press.







  \NAME{M.\ Pfender} 2008 RCF1: Theories of PR Maps and Partial PR Maps.
  pdf file. TU Berlin.

  \NAME{M.\ Pfender} 2008 RCFX: Universal Objects and Theory Embedding.
  pdf file. TU Berlin.

  \NAME{M.\ Pfender, M.\ Kr\"oplin, D.\ Pape} 1994: Primitive
  Recursion, Equality, and a Universal Set. 
  \emph{Math.\ Struct.\ in Comp.\ Sc.\ } \textbf{4,} 295-313.
  
  

  \NAME{W.\ Rautenberg} 1995/2006: \emph{A Concise Introduction to 
  Mathematical Logic.} Universitext Springer 2006.


  \NAME{R.\ Reiter} 1980: Mengentheoretische Konstruktionen in arithmetischen
  Universen. Diploma Thesis. TU Berlin.

  
  \NAME{L.\  Rom\`an} 1989: Cartesian categories with natural numbers object.
  \emph{J.\  Pure and Appl.\  Alg.} \textbf{58,} 267-278.
  



  \NAME{W.\ W.\ Tait} 1996: Frege versus Cantor and Dedekind: on the concept
  of number. Frege, Russell, Wittgenstein: \emph{Essays in Early Analytic
  Philosophy (in honor of Leonhard Linsky)} (ed. W.\ W.\ Tait). Lasalle:
  Open Court Press (1996): 213-248. Reprinted in \emph{Frege: Importance 
  and Legacy} (ed. M.\ Schirn). Berlin: Walter de Gruyter (1996): 70-113.


  \NAME{A.\ Tarski, S.\ Givant} 1987: \emph{A formalization of set theory 
  without variables}. AMS Coll.\ Publ.\ vol.\ 41.

  \NAME{M.\ Tierney} 1973: Axiomatic Sheaf Theory. 
  \emph{C.\ I.\ M.\ E.\ Conf.\ on Categories and Commutative Algebra,
  Varenna,} pp.\ 249-326. Edizione Cremonese. Roma.
  (Quoted in \NAME{Wraith} 1973.)

  \NAME{G.\ Wraith} 1973: Lectures on Elementary Topoi. In
  \emph{Model Theory and Topoi,} LN in Math.\ \textbf{445,} 114-206.

\bigskip

\bigskip


  \noindent Address of the author: \\
  \NAME{M. Pfender}                       \hfill D-10623 Berlin \\
  Institut f\"ur Mathematik                              \\
  Technische Universit\"at Berlin         \hfill pfender@math.TU-Berlin.DE\\
  

\vfill

\end{document}